\documentclass[12pt,leqno]{amsart}
\usepackage{amssymb,amsthm,verbatim,enumerate,ifthen,xcolor,hyperref}
\usepackage[mathscr]{eucal}
\usepackage[utf8]{inputenc}
\usepackage[T1]{fontenc}
\usepackage{graphics}
\usepackage{graphicx}
\def\N{\mathbb{N}}
\def\R{\mathbb{R}}
\def\Q{\mathbb{Q}}
\def\Z{\mathbb{Z}}

\def\sup{\mathop{\mbox{\rm sup}}}
\def\dist{\mathop{\mbox{\rm dist}}}

\def\supp{\mathop{\mbox{\rm supp}}}

\def\T{\mathscr{T}}

\def\B{\mathscr{B}}
\def\epsilon{\varepsilon}
\def\eps{\varepsilon}
\def\de{\delta}

\newtheorem{theorem}{Theorem}
\newtheorem*{theorem*}{Theorem}
\def\Thm#1#2{\ifthenelse{\equal{#1}{*}}{\begin{theorem*}#2\end{theorem*}}
  {\begin{theorem}\label{T#1}#2\end{theorem}}}
\newtheorem{Atheorem}{Theorem}

\def\thm#1{Theorem~\ref{T#1}}

\newtheorem{proposition}[theorem]{Proposition}
\newtheorem*{proposition*}{Proposition}
\def\Prp#1#2{\ifthenelse{\equal{#1}{*}}{\begin{proposition*}#2\end{proposition*}}
             {\begin{proposition}\label{P#1}#2\end{proposition}}}
\def\prp#1{Proposition~\ref{P#1}}
\newtheorem{corollary}[theorem]{Corollary}
\newtheorem*{corollary*}{Corollary}
\def\Cor#1#2{\ifthenelse{\equal{#1}{*}}{\begin{corollary*}#2\end{corollary*}}
             {\begin{corollary}\label{C#1}#2\end{corollary}}}

\newtheorem{lemma}[theorem]{Lemma}
\newtheorem*{lemma*}{Lemma}
\def\Lem#1#2{\ifthenelse{\equal{#1}{*}}{\begin{lemma*}#2\end{lemma*}}
             {\begin{lemma}\label{L#1}#2\end{lemma}}}
\def\lem#1{Lemma~\ref{L#1}}

\newtheorem{example}[theorem]{Example}
\newtheorem*{example*}{Example}
\long\def\Exa#1#2{\ifthenelse{\equal{#1}{*}}{\begin{example*}\rm #2\end{example*}}
            {\begin{example}\label{Ex#1}\rm #2\end{example}}}

\newtheorem{problem}[subsection]{Problem}

\theoremstyle{definition}
\newtheorem{definition}[theorem]{Definition}
\newtheorem*{definition*}{Definition}
\def\Defi#1#2{\ifthenelse{\equal{#1}{*}}{\begin{definition*}#2\end{definition*}}
      {\begin{definition}\label{D#1}#2\end{definition}}}

\newtheorem{remark}[theorem]{Remark}
\newtheorem*{remark*}{Remark}
\def\Rem#1#2{\ifthenelse{\equal{#1}{*}}{\begin{remark*}#2\end{remark*}}
             {\begin{remark}\label{R#1}#2\end{remark}}}
\def\rem#1{Remark~\ref{R#1}}

\def\eq#1{{\rm(\ref{E#1})}}
\def\Eq#1#2{\ifthenelse{\equal{#1}{*}}
  {\begin{equation*}\begin{aligned}#2\end{aligned}\end{equation*}}
  {\begin{equation}\begin{aligned}\label{E#1}#2\end{aligned}\end{equation}}}

\begin{document}
\begin{flushright}
Publ. Math. Debrecen \textbf{80}(1-2) (2012), 107–126. \\
\href{http://dx.doi.org/110.5486/PMD.2012.4930}{110.5486/PMD.2012.4930} \\[1cm]
\end{flushright}

\title[]{On $\varphi$-convexity}
\author[J.\ Makó]{Judit Makó}
\author[Zs. Páles]{Zsolt Páles}
\address{Institute of Mathematics, University of Debrecen,
H-4010 Debrecen, Pf.\ 12, Hungary}
\email{\{makoj,pales\}@science.unideb.hu}


\thanks{This research has been supported by the Hungarian Scientific Research Fund (OTKA) Grant
NK81402 and by the TÁMOP 4.2.1./B-09/1/KONV-2010-0007 project implemented through the New
Hungary Development Plan co-financed by the European Social Fund, and the European Regional
Development Fund.}

\begin{abstract}
In this paper, approximate convexity and approximate midconvexity properties,
called $\varphi$-convexity and $\varphi$-midconvexity, of real valued function are investigated.
Various characterizations of $\varphi$-convex and $\varphi$-midconvex functions are obtained.
Furthermore, the relationship between $\varphi$-midconvexity and $\varphi$-convexity is established.
\end{abstract}

\maketitle

\section{Introduction}

The stability theory of functional inequalities started with the paper
\cite{HyeUla52} of Hyers and Ulam who introduced the notion of $\eps$-convex
function: If $D$ is a convex subset of a real linear space $X$ and
$\eps$ is a nonnegative number, then a function $f:D\to\R$ is called
{\it $\eps$-convex} if
\Eq{0}{
  f(tx+(1-t)y)\le tf(x)+(1-t)f(y) + \eps
}
for all $x,y\in D$, $t\in[0,1]$.  The basic result obtained by Hyers and
Ulam states that if the underlying space $X$ is of finite dimension then
$f$ can be written as $f=g+h$, where $g$ is a convex function and $h$ is
a bounded function whose supremum norm is not larger than $k_n\eps$,
where the positive constant $k_n$ depends only on the dimension $n$ of
the underlying space $X$.  Hyers and Ulam proved that
$k_n\le(n(n+3))/(4(n+1))$. Green \cite{Gre52}, Cholewa \cite{Cho84a}
obtained much better estimations of $k_n$ showing that asymptotically
$k_n$ is not bigger than $(\log_2(n))/2$. Laczkovich \cite{Lac99}
compared this constant to several other dimension-depending stability
constants and proved that it is not less than $(\log_2(n/2))/4$. This
result shows that there is no analogous stability results for infinite dimensional
spaces $X$. A counterexample in this direction was earlier constructed
by Casini and Papini \cite{CasPap93}. The stability aspects of
$\eps$-convexity are discussed by Ger \cite{Ger94e}. An overview of results on
$\delta$-convexity can be found in the book of Hyers, Isac, and
Rassias \cite{HyeIsaRas98a}.

If $t=1/2$ and \eq{0} holds for all $x,y\in D$, then $f$ is called an
{\it $\eps$-Jensen-convex function}.  There is no analogous decomposition
for $\eps$-Jensen-convex functions by the counterexample given by Cholewa
\cite{Cho84a}.  However, one can get Bernstein-Doetsch type regularity
theorems which show that $\eps$-Jensen-convexity and local upper
boundedness imply $2\eps$-convexity.  This result is due to Bernstein
and Doetsch \cite{BerDoe15} for $\eps=0$, and to Ng and Nikodem \cite{NgNik93}
in the case $\eps\ge0$.  For some recent extensions of these results
to more general convexity concepts, see \cite{Pal00b}. For locally
upper bounded $\eps$-Jensen-convex functions one can obtain the existence
of an analogous stability constant $j_n$ (defined similarly as $k_n$
above). The sharp value of this stability constant has recently been
found by Dilworth, Howard, and Roberts \cite{DilHowRob99} who have shown that
$$
  j_n=\frac{1}{2}\Bigl([\log_2(n)]+1+\frac{n}{2^{[\log_2(n)]}}\Bigr)
     \le 1+\frac{1}{2}\log_2(n)
$$
is the best possible value for $j_n$. (Here $[\cdot]$ denotes the
integer-part function). The connection between $\eps$-Jensen-convexity
and $\eps$-$\Q$-convexity has been investigated by Mrowiec \cite{Mro01}.

If $D\subset\R$ and \eq{0} is supposed to be valid for all $x,y\in D$
except a set of 2-dimensional Lebesgue measure zero then one can speak
about {\it almost $\eps$-convexity}.  Results in this direction are due
to Kuczma \cite{Kuc70a} (the case $\eps=0$) and Ger \cite{Ger88c} (the
case $\eps\ge0$).

In a recent paper \cite{Pal03a}, the second author introduced a more
general notion than $\eps$-convexity. Let $\eps$
and $\de$ be nonnegative constants. A function $f:D \to\R$
is called $(\eps,\delta)$-convex, if
$$
  f\left(t x+(1-t) y \right)
    \leq t f(x) + (1-t) f(y) + \delta + \eps t (1-t) \|x-y\|
$$
for every $x,y\in D$ and $t \in[0,1]$. The main results of the paper
\cite{Pal03a} obtain a complete characterization of $(\eps,\delta)$-convexity
if $D\subseteq\R$ is an open real interval
by showing that these functions are of the form $f=g+h+\ell$, where $g$ is convex,
$h$ is bounded with $\|h\|\leq \de/2$ and $\ell$ is Lipschitzian with Lipschitz modulus
Lip$(\ell)\leq\eps$.

In the papers \cite{HazPal04}, \cite{HazPal05}, the notion of $(\eps,p)$-convexity and
$(\eps,p)$-midconvexity were introduced: If $\eps,p\geq0$ and $t\in[0,1]$, then
a function $f:D \to\R$ is called \textit{$(\eps,p,t)$-convex}, if
$$
  f\left(t x+(1-t) y \right)
    \leq t f(x) + (1-t) f(y) + \eps (t (1-t) \|x-y\|)^p
$$
for every $x,y\in D$. If the above property holds for $t=1/2$ and for all $t\in[0,1]$,
then we speak about \textit{$(\eps,p)$-midconvexity} and \textit{$(\eps,p)$-convexity}, 
respectively. The main result in \cite{HazPal05} shows that, for locally upper bounded 
functions, $(\eps,p)$-midconvexity implies $(c\eps,p)$-convexity for some constant $c$.

Another, but related, notion of approximate convexity, the concept of so-called 
paraconvexity was introduced by Rolewicz \cite{Rol79b,Rol79a,Rol05b} in the late 70s.
It also turned out that Takagi-like functions appear naturally in the
investigation of approximate convexity, see, for example, Boros \cite{Bor08},
Házy \cite{Haz07a,Haz07b}, Házy and Páles \cite{HazPal04,HazPal05,HazPal09},
Makó and Páles \cite{MakPal11b,MakPal11a}, 
Mrowiec, Tabor and Tabor \cite{MroTabTab08}, Tabor and Tabor \cite{TabTab09b,TabTab09a},
Tabor, Tabor, and Żołdak \cite{TabTabZol10a,TabTabZol10b}.

The aim of this paper is to offer a unified framework for most of the mentioned
approximate convexity notions by introducing the notions of $\varphi$-convexity and
$\varphi$-midconvexity and to extend the previously known results to this more general setting.
We also introduce the relevant Takagi type functions which appear naturally in the
description of the connection of $\varphi$-convexity and $\varphi$-midconvexity.

\section{$\varphi$-convexity and $\varphi$-midconvexity}

Throughout the paper $\R$, $\R_+$, and $\N$ denote the sets of real, nonnegative
real, and natural numbers, respectively.
Assume that $D$ is a nonempty convex subset of a real normed space $X$
and denote $D^+:=\{\|x-y\|:x,y\in D\}$. Let $\varphi:D^+\to \R_+$ be a given function.

\Defi{1}{A function $f:D\to \R$ is called \textit{$\varphi$-convex} on $D$, if
\Eq{1a}
{
f(tx+(1-t)y)\leq tf(x)+(1-t)f(y)+t\varphi\big((1-t)\|x-y\|\big)+(1-t)\varphi\big(t\|x-y\|\big)
}
holds for all $t\in [0,1]$ and for all $x,y \in D.$ If \eq{1a} holds for $t=1/2$,
i.e., if, for all $x,y \in D,$
\Eq{1b}
{
f\bigg(\frac{x+y}{2}\bigg)\leq \frac{f(x)+f(y)}{2}+\varphi\bigg(\Big\|\frac{x-y}{2}\Big\|\bigg),
}
then we say that $f$ is \textit{$\varphi$-midconvex}.
}

In the case $\varphi\equiv0$, the meaning of inequalities \eq{1a} and \eq{1b} is the convexity
and midconvexity (Jensen-convexity) of $f$, respectively.

An important particular case occurs when $\varphi:D^+\to\R_+$ is of the form
$\varphi(x):=\epsilon x^p,$ where $p,\epsilon\geq0$ are arbitrary constants.
Then the function $f$ is called \textit{$(\epsilon,p)$-convex} and
\textit{$(\epsilon,p)$-midconvex} on $D$, respectively (cf.\ \cite{Pal03a}).

The next results describe the structure of $\varphi$-convex functions and
$\varphi$-midconvex functions.

\Prp{2}{{\color{white}.}
\begin{enumerate}[(i)]
    \item If, for $j=1,\dots,n$, $\varphi_j:D^+\to\R_+$, the function $f_j:D\to \R$ is
    $\varphi_j$-convex and $c_j$ is a nonnegative number,
    then $c_1f_1+\cdots+c_nf_n$ is $(c_1\varphi_1+\cdots+c_n\varphi_n)$-convex.
    In particular, the set of $\varphi$-convex functions on $D$ is convex.
    \item Let $\{f_{\gamma}:D\to\R\mid\gamma\in \Gamma \}$ be a family of
    $\varphi$-convex functions. Assume, for all $x\in D$, that
    $f(x):=\sup_{\gamma\in \Gamma}f_{\gamma}(x)<+\infty$. Then $f$ is $\varphi$-convex.
    \item Let $\{f_{\gamma}:D\to\R\mid\gamma\in \Gamma \}$ be a downward directed
    family of $\varphi$-convex functions in the following sense:
    for all $\gamma_1,\gamma_2\in\Gamma$ and $x_1,x_2\in D$, there exists $\gamma\in\Gamma$
    such that $f_\gamma(x_i)\leq f_{\gamma_i}(x_i)$ for $i=1,2$.
    Assume, for all $x\in D$, that $f(x):=\inf_{\gamma\in \Gamma}f_{\gamma}(x)>-\infty$.
    Then $f$ is $\varphi$-convex.
\end{enumerate}}

\begin{proof}
(i) is easy to prove.

(ii) Let $x,y\in D$ and $t\in[0,1]$. For all $\gamma\in \Gamma,$ we have
\Eq{*}
{
f_{\gamma}(tx+(1-t)y)
   &\leq tf_{\gamma}(x)+(1-t)f_{\gamma}(y)+ t\varphi\big((1-t)\|x-y\|\big)+(1-t)\varphi\big(t\|x-y\|\big)\\
   &\leq tf(x)+(1-t)f(y)+ t\varphi\big((1-t)\|x-y\|\big)+(1-t)\varphi\big(t\|x-y\|\big).}
Thus,
\Eq{*}
{
f(tx+(1-t)y)&=\sup_{\gamma\in \Gamma}f_{\gamma}(tx+(1-t)y)\\
            &\leq tf(x)+(1-t)f(y)+t\varphi\big((1-t)\|x-y\|\big)+(1-t)\varphi\big(t\|x-y\|\big).
}
Hence $f$ is $\varphi$-convex.

(iii) Let $x,y \in D$ and $t\in [0,1].$ Let $\delta>0$ be arbitrary.
Then $f(x)<f(x)+\delta$ and $f(y)<f(y)+\delta.$
Thus there exist $\gamma_1,\gamma_2,$ such that $f_{\gamma_1}(x)<f(x)+\delta$ and $f_{\gamma_2}(y)<f(y)+\delta.$
By the conditions of the proposition, there exists $\gamma\in \Gamma,$ such that
\Eq{*}
{
f_{\gamma}(x)&\leq f_{\gamma_1}(x)<f(x)+\delta,\\
f_{\gamma}(y)&\leq f_{\gamma_2}(y)<f(y)+\delta.
}
Then we get
\Eq{*}
{
f(tx+(1-t)y)&\leq f_{\gamma}(tx+(1-t)y)\\
            &\leq tf_{\gamma}(x)+(1-t)f_{\gamma}(y)+t\varphi\big((1-t)\|x-y\|\big)+(1-t)\varphi\big(t\|x-y\|\big)\\
            &\leq tf(x)+(1-t)f(y)+\delta+t\varphi\big((1-t)\|x-y\|\big)+(1-t)\varphi\big(t\|x-y\|\big).
}
This proves that $f$ is $\varphi$-convex.
\end{proof}

The following statements concern midconvex funtions, they are analogous to those
of \prp{2}.

\Prp{3}{{\color{white}.}
\begin{enumerate}[(i)]
    \item If, for $j=1,\dots,n$, $\varphi_j:D^+\to\R_+$, the function $f_j:D\to \R$ is
    $\varphi_j$-midconvex and $c_j$ is a nonnegative number,
    then $c_1f_1+\cdots+c_nf_n$ is $(c_1\varphi_1+\cdots+c_n\varphi_n)$-midconvex.
    In particular, the set of $\varphi$-midconvex functions on $D$ is convex.
    \item Let $\{f_{\gamma}:D\to\R\mid\gamma\in \Gamma \}$ be a family of
    $\varphi$-midconvex functions. Assume, for all $x\in D$, that
    $f(x):=\sup_{\gamma\in \Gamma}f_{\gamma}(x)<+\infty$. Then $f$ is $\varphi$-midconvex.
    \item Let $\{f_{\gamma}:D\to\R\mid\gamma\in \Gamma \}$ be a downward directed
    family of $\varphi$-midconvex functions in the following sense:
    for all $\gamma_1,\gamma_2\in\Gamma$ and $x_1,x_2\in D$, there exists $\gamma\in\Gamma$
    such that $f_\gamma(x_i)\leq f_{\gamma_i}(x_i)$ for $i=1,2$.
    Assume, for all $x\in D$, that $f(x):=\inf_{\gamma\in \Gamma}f_{\gamma}(x)>-\infty$.
    Then $f$ is $\varphi$-midconvex.
\end{enumerate}}

\Defi{H}{A function $f:D\to \R$ is said to be of \textit{$\varphi$-Hölder class} on $D$
or briefly $f$ is called \textit{$\varphi$-Hölder} on $D$ if there exists a
nonnegative constant $H$ such that, for all $x,y\in D,$
\Eq{H}
{
|f(x)-f(y)|\leq H\varphi(\|x-y\|).
}
The smallest constant $H$ such that \eq{H} holds is said to be the
\textit{$\varphi$-Hölder modulus} of $f$ and is denoted by ${H}_{\varphi}(f)$.
}

A relationship between the $\varphi$-Hölder property and $\varphi$-convexity is
obtained in the following result.

\Prp{H}{Let $f:D\to \R$ be of $\varphi$-Hölder class on $D$.
Then $f$ is $\big({H}_{\varphi}(f)\cdot\varphi\big)$-convex on $D$.}
\begin{proof}
Let $x,y\in D$ and let $t\in [0,1].$ Then
\Eq{*}
{
f(tx+(1-t)y)&-tf(x)-(1-t)f(y)\\
   &=t\big(f(tx+(1-t)y)-f(x)\big)+(1-t)\big(f(tx+(1-t)y)-f(y)\big)\\
   &\leq t{H}_{\varphi}(f)\varphi(\|tx+(1-t)y-x\|)+ (1-t){H}_{\varphi}(f)\varphi(\|tx+(1-t)y-y\|),
}
which is equivalent to the $\big({H}_{\varphi}(f)\cdot\varphi\big)$-convexity of $f.$
\end{proof}

For functions $\varphi:D^+\to \R$, we introduce the following subadditivity-type property:

\Defi{IS}{
We say that $\varphi$ is \textit{increasingly subadditive on $D^+$} if,
for all $u,v,w\in D^+$ with $u\leq v+w$,
\Eq{uvw}{
  \varphi(u)\leq\varphi(v)+\varphi(w)
}
holds.} Clearly, if $\varphi:\R_+\to \R$ is nondecreasing and subadditive then it is also
increasingly subadditive on $\R^+$.

\Prp{5}{Assume that $\varphi:D^+\to \R$ is increasingly subadditive. Then, for all $z\in D$,
the map $x \mapsto -\varphi(\|x-z\|)$ is of $\varphi$-Hölder class on $D$ with
$\varphi$-Hölder modulus $1$, and therefore, it is also $\varphi$-convex on $D$.}

\begin{proof} Let $z\in D$ be fixed.
To prove the $\varphi$-Hölder property of the map $x \mapsto -\varphi(\|x-z\|)$, let $x,y\in D.$
Then $u=\|x-z\|$, $v=\|x-y\|$, and $w=\|y-z\|$ are elements of $D^+$ such that \eq{uvw} holds.
Therefore, by the increasing subadditivity, we get
\Eq{*}
{
\varphi(\|x-z\|)-\varphi(\|y-z\|)
\leq \varphi(\|x-y\|)+\varphi(\|y-z\|)-\varphi(\|y-z\|)=\varphi(\|x-y\|).
}
Interchanging $x$ and $y$, we also have $\varphi(\|y-z\|)-\varphi(\|x-z\|)\leq \varphi(\|y-x\|).$
These two inequalities imply
\Eq{*}
{
\big|\varphi(\|x-z\|)-\varphi(\|y-z\|)\big|\leq \varphi(\|x-y\|),
}
which means that the map $x\mapsto -\varphi(\|x-z\|)$ is $\varphi$-Hölder on $D$
with $\varphi$-Hölder modulus $1$.
\end{proof}

The next lemma is well known, for completeness we provide its short proof.

\Lem{4}{Let $0\leq p\leq 1$ be an arbitrary constant. Then the map $x\mapsto x^p$ is
subadditive and nondecreasing on $\R_+$ and hence it is also increasingly subadditive on $\R^+$.}

\begin{proof} For $s\in ]0,1[$, we have $s\leq s^p.$ Hence
\Eq{*}{
    1=s+(1-s)\leq s^p+(1-s)^p.
}
If $x,y\in\R_+$ then, with $s:=\frac{x}{x+y}\in ]0,1[$, we get
\Eq{*}{
  1\leq \Big(\frac{x}{x+y}\Big)^p+\Big(\frac{y}{x+y}\Big)^p,
}
which shows the subadditivity of the function $x\mapsto x^p$.
\end{proof}

\Defi{5}{
Let $0<p\leq 1$ be an arbitrary constant.
For all $t\in D^+$ let $\varphi(t):=t^p,$ then if $f:D\to \R$ is a $\varphi$-Hölder function,
then it is called (classical) \textit{$p$-Hölder functions}.
In this case the $\varphi$-Hölder modulus is called \textit{$p$-Hölder modulus} of $f$
and it is denoted by ${H}_{p}(f)$.
}

The next corollary gives a relationship between the $p$-Hölder
functions and the $p$-convex functions.

\Cor{4}{
Let $0<p\leq 1$ be an arbitrary constant and $z\in X$. Then $x \mapsto -\|x-z\|^p$ is
of $p$-Hölder class on $X$ with the $p$-Hölder modulus $1$, and therefore,
it is $(1,p)$-convex on $X$.
}

The subsequent theorem, which is one of the main results of this paper,
offers equivalent conditions for $\varphi$-convexity. It generalizes the result
of \cite[Thm. 1]{Pal03a}.

\Thm{2}{Let $D$ be an open real interval and $f:D\to \R.$ Then the following conditions are equivalent.
\begin{enumerate}[(i)]
 \item $f$ is $\varphi$-convex on $D.$
 \item For $x,u,y\in D$ with $x<u<y,$
   \Eq{2b}
   {
   \frac{f(u)-f(x)-\varphi(u-x)}{u-x}\leq \frac{f(y)-f(u)+\varphi(y-u)}{y-u}.
   }
 \item There exists a function $a:D\to \R$ such that, for $x,u\in D,$
   \Eq{2c}
   {
   f(x)-f(u)\geq a(u)(x-u)-\varphi(|x-u|).
   }
\end{enumerate}}

\begin{proof}
$(\mbox{i}) \Rightarrow (\mbox{ii})$ Assume that $f:D\to \R$ is
$\varphi$-convex and let $x<u<y$ be arbitrary elements of $D.$
 Choose $t\in [0,1]$ such that $u=tx+(1-t)y,$ that is let $t:=\frac{y-u}{y-x}.$
 Then, applying the $\varphi$-convexity of $f$, we get
\Eq{*}
{
f(u)\leq \frac{y-u}{y-x}f(x)+\frac{u-x}{y-x}f(y)+\frac{y-u}{y-x}\varphi\Big(\frac{u-x}{y-x}(y-x)\Big)
+\frac{u-x}{y-x}\varphi\Big(\frac{y-u}{y-x}(y-x)\Big),
}
which is equivalent to
\Eq{*}
{
(y-u)\big(f(u)-f(x)-\varphi(u-x)\big)\leq (u-x)\big(f(y)-f(u)+\varphi(y-u)\big).
}
Dividing by $(y-u)(u-x)>0$, we arrive at \eq{2b}.

$(\mbox{ii}) \Rightarrow (\mbox{iii})$ Assume that $(\mbox{ii})$
holds and, for $u\in D$, define \Eq{*} { a(u):= \inf_{y\in D,u<y}
\frac{f(y)-f(u)+\varphi(y-u)}{y-u}. } Then in view of
$(\mbox{ii}),$ we get \Eq{2d} {
\frac{f(x)-f(u)-\varphi(u-x)}{x-u}\leq a(u)\leq
\frac{f(y)-f(u)+\varphi(y-u)}{y-u}, } for all $x<u<y$ in $D.$ The
left-hand side inequality in $\eq{2d}$ yields $\eq{2c}$ in the
case $x<u,$ and analogously, the right-hand side inequality (with
the substitution $y:=x$) reduces to $\eq{2c}$ in the case $x>u.$
The case $x=u$ is obvious.

$(\mbox{iii}) \Rightarrow (\mbox{i})$
Let $x,y\in D,$ $t\in [0,1]$, and set $u:=tx+(1-t)y$. Then, by $(\mbox{iii})$, we have
\Eq{*}
{
f(x)-f(u)&\geq a(u)(x-u)-\varphi(|x-u|),\\
f(y)-f(u)&\geq a(u)(y-u)-\varphi(|y-u|).
}
Multiplying the first inequality by $t$ and the second inequality by $1-t$
and adding up the inequalities so obtained, we get $\eq{1a}.$
\end{proof}

\Rem{12}{In the vector variable setting (i.e., when $D$ is an open convex subset of a normed space $X$),
instead of condition (iii), the following analogous property can be formulated:
\begin{enumerate}
 \item[(iii)$^*$] There exists a function $a:D\to X^*$ such that, for $x,u\in D,$
   \Eq{*}
   {
   f(x)-f(u)\geq a(u)(x-u)-\varphi(\|x-u\|).
   }
\end{enumerate}
One can easily see (by the same argument as above) that (iii)$^*$ implies (i), that is,
the $\varphi$-convexity of $f$. The validity of the reversed implication is an open problem.}

The next theorem gives another characterization of $\varphi$-convex functions,
if $\varphi$ is increasingly subadditive.

\Thm{6}{Let $D$ be an open real interval and let $\varphi:D^+\to \R_+$ be increasingly subadditive.
Then a function $f:D\to\R$ is $\varphi$-convex if and only if there exist
two functions $a:D\to \R$ and $b:D\to \R$ such that
\Eq{6a}
{
f(x)=\sup_{u\in D}\big(a(u)x+b(u)-\varphi(|x-u|)\big),
}
for all $x\in D.$}

\begin{proof} Assume that $f$ is $\varphi$-convex. By \thm{2}, there exists
a function $a:D\to \R$ such that
\Eq{*}
{
f(x)\geq f(u)+a(u)(x-u)-\varphi(|x-u|),
}
for all $u,x\in D.$
Define $b(u):=f(u)-a(u)u,$ for $u\in D.$ Thus, for $u,x\in D$,
\Eq{*}
{
f(x)\geq a(u)x+b(u)-\varphi(|x-u|)
}
and we have equality for $u=x$. Therefore, \eq{6a} holds.

Conversely, assume that \eq{6a} is valid for $x\in D$.
By \prp{5}, for fixed $u\in D$, the mapping $x\mapsto-\varphi(|x-u|)$ is
$\varphi$-convex. The map $x\mapsto a(u)x+b(u)$ is affine, and hence
the function $f_u:D\to\R$ defined by $f_u(x):=a(u)x+b(u)-\varphi(|x-u|)$ is
$\varphi$-convex for all fixed $u\in D$.
Now applying (ii) of \prp{2}, we obtain that $f$ is $\varphi$-convex.
\end{proof}

\Rem{13}{In the vector variable setting (i.e., when $D$ is an open convex subset of a normed space $X$),
the following implication can be formulated: If $\varphi:D^+\to \R_+$ is increasingly subadditive
and there exist two function $a:D\to X^*$ and $b:D\to \R$ such that, for $x\in D,$
\Eq{*}
{
f(x)=\sup_{u\in D}\big(a(u)(x)+b(u)-\varphi(\|x-u\|)\big),
}
then $f$ is $\varphi$-convex. The validity of the reversed implication is an open problem.}

\Cor{11}
{Let $D$ be an open real interval and
let $0< p\leq1$ and $\epsilon\geq0$ be arbitrary constants. Then a function $f:D\to\R$ is
$(\epsilon,p)$-convex if and only if there exist two functions
$a:D\to \R$ and $b:D\to \R$ such that
\Eq{*}
{
f(x)=\sup_{u\in D}\big(a(u)x+b(u)-\epsilon|x-u|^p\big),
}
for all $x\in D.$
}

The subsequent theorem offers a sufficient condition for the $\varphi$-midconvexity.
The result is analogous to the implication (iii)$\Rightarrow$(i) of \thm{2}.
Unfortunately, we were not able to obtain the necessity of this condition, i.e.,
the reversed implication.

\Thm{2M}{Let $f:D\to \R$ and assume that, for all $u\in D$,
there exists an additive function $A_u: X\to X$ such that
\Eq{2Ma}
{
f(x)-f(u)\geq A_u(x-u)-\varphi(\|x-u\|) \qquad(x\in D).
}
Then, $f$ is $\varphi$-midconvex.}

\begin{proof}
Let $x,y\in D$ and set $u:=\frac{x+y}{2}$. Then, by \eq{2Ma}, we have
\Eq{*}
{
f(x)-f(u)&\geq A_u(x-u)-\varphi(\|x-u\|)=A_u\bigg(\frac{x-y}{2}\bigg)-\varphi \bigg(\Big\|\frac{x-y}{2}\Big\| \bigg),\\
f(y)-f(u)&\geq A_u(y-u)-\varphi(\|y-u\|)=A_u\bigg(\frac{y-x}{2}\bigg)-\varphi \bigg(\Big\|\frac{y-x}{2}\Big\| \bigg).
}
Adding up the inequalities and multiplying the inequality so obtained by $\frac{1}{2}$, we get $\eq{1b}.$
\end{proof}

The following result is analogous to \thm{6}, however it offers only a sufficient
condition for $\varphi$-midconvexity.

\Thm{6M}{Let $\varphi:D^+\to \R_+$ be increasingly subadditive and let $f:D\to \R$.
Assume that, for all $u\in D$, there exists an additive function $A_u:X\to X$ and
there exists a function $b:D\to\R$ such that
\Eq{6M}
{
f(x)=\sup_{u\in D}\big(A_u(x)+b(u)-\varphi(\|x-u\|)\big),
}
for all $x\in D.$ Then $f$ is $\varphi$-midconvex.}

\begin{proof} Assume that \eq{6M} is valid for $x\in D$.
By \prp{5}, for fixed $u\in D$, the mapping $x\mapsto-\varphi(\|x-u\|)$ is
$\varphi$-convex, so it is $\varphi$-midconvex. The map $x\mapsto A_u(x)+b(u)$ is affine,
 and hence the function $f_u:D\to\R$ defined by $f_u(x):=A_u(x)+b(u)-\varphi(\|x-u\|)$ is
$\varphi$-midconvex for all fixed $u\in D$.
Now applying (ii) of \prp{3}, we obtain that $f$ is $\varphi$-midconvex.
\end{proof}

Henceforth we search for relations between the local upper-bounded
$\varphi$-midconvex functions and $\varphi$-convex functions
with the help of the results from the papers \cite{Haz05a} and \cite{HazPal05}
by Házy and Páles.

Define the function $d_\Z:\R\to\R_+$ by
\Eq{*}{
  d_\Z(t)=\dist(t,\Z):=\min\{|t-k|:k\in\Z\}.
}
It is immediate to see that $d_\Z$ is 1-periodic and symmetric with respect to
$t=1/2$, i.e., $d_\Z(t)=d_\Z(1-t)$ holds for all $t\in\R$.
For a fixed $\varphi:\frac12 D^+\to\R_+$, we introduce the Takagi type
function $\T_\varphi:\R\times D^+\to\R_+$ by
\Eq{TT}
{
\T_\varphi(t,u):=\sum_{n=0}^{\infty}\frac{\varphi\big(d_\Z(2^{n}t)u\big)}{2^{n}}
\qquad((t,u)\in \R\times D^+).
}
Applying the estimate $0\leq d_\Z\leq\frac12$,
one can easily see that $\T_\varphi(t,u)\leq2\varphi\big(\frac{u}{2}\big)$ for $u\in D^+$
whenever $\varphi$ is nondecreasing.

For $p\geq0$, we also define the Takagi type function $T_p:\R\to\R_+$ by
\Eq{Tp}
{
 T_p(t):=\sum_{n=0}^{\infty}\frac{\big(d_\Z(2^{n}t)\big)^p}{2^{n}}
\qquad(t\in\R).
}
In the case when $\varphi$ is of the form $\varphi(t)=\epsilon|t|^p$ for some
constants $\epsilon\geq0$ and $p\geq0$, the following identity holds:
\Eq{*}{
   \T_\varphi(t,u)=\epsilon T_p(t)u^p \qquad((t,u)\in \R\times D^+).
}
Observe that $\T_\varphi$ and $T_p$ are also 1-periodic and symmetric with respect to
$t=1/2$ in their first variables.

In order to obtain lower and upper estimates for the functions $\T_\varphi$ and
$T_p$ defined above, we need to recall de Rham's classical theorem \cite{Rha57}.
By $\B(\R,\R)$ we denote the space of bounded functions $f:\R\to\R$ equipped
with the supremum norm.

\Thm{TT}{Let $\psi\in \B(\R,\R), a,b\in \R,$ $|a|<1.$ Let
$F_\psi:\B(\R,\R)\to \B(\R,\R)$ be an operator defined as follows
\Eq{*}
{
\big(F_\psi f\big)(t)
  :=af(bt)+\psi(t) \qquad \mbox{for}\quad f\in \B(\R,\R),\,\,t\in \R.
}
Then
\begin{enumerate}[(i)]
\item $F_\psi$ is a contraction on $\B(\R,\R)$ with a unique fixed point
$f_\psi$ which is given by the formula
\Eq{*}{
  f_\psi(t)=\sum_{n=0}^{\infty}a^n\psi(b^nt) \qquad (t\in \R);
}
\item if $a\geq 0$ and the functions $g,h\in\B(\R,\R)$ satisfy
the inequalities $g\leq F_\psi g$ and $F_\psi h\leq h$, then $g\leq f_\psi\leq h.$
\end{enumerate}}

\Rem{1}{In view of the first assertion of this theorem, observe that the functions
$\T_\varphi(\cdot,u)$ and $T_p$
defined in \eq{TT} and \eq{Tp} are the fixed points of the operator:
\Eq{FF}{
   \big(F_\psi f\big)(t):=\frac12 f(2t)+\psi(t)
    \qquad \mbox{for}\quad f\in \B(\R,\R),\,\,t\in\R
}
where $\psi\in\B(\R,\R)$ is given by $\psi(t):=\varphi\big(d_\Z(t)u\big)$
and $\psi(t):=\big(d_\Z(t)\big)^p$, respectively.}

In the results below, we establish upper
and lower bounds for $\T_\varphi$ in terms of the function $\tau_\varphi:\R\times D^+\to\R$ defined by
\Eq{*}
{
\tau_{\varphi}(t,u):= d_\Z(t)\varphi\big((1-d_\Z(t))u\big)
         +(1-d_\Z(t))\varphi\big(d_\Z(t)u\big) \qquad ((t,u)\in \R\times D^+).
}
Observe that, for $t\in[0,1]$, we have
\Eq{*}
{
\tau_{\varphi}(t,u):= t\varphi\big((1-t)u\big)
         +(1-t)\varphi\big(tu\big) \qquad (u\in D^+),
}
which is exactly the error term related to $\varphi$-convexity.

\Prp{7}{Let $\varphi:D^+\to \R_+$ be subadditive. Then, for all $(t,u)\in\R\times D^+,$
\Eq{7a}
{
\tau_\varphi(t,u)\leq \T_\varphi(t,u).
}}

\begin{proof}
Let $u\in D^+$ be arbitrarily fixed. By the 1-periodicity and symmetry with respect to
the point $t=1/2$, it suffices to show that \eq{7a} holds for all $t\in\big[0,\frac12\big]$.
If $t=0$ then \eq{7a} is obvious.
Now assume that $0<t\leq \frac12.$
Then there exists a unique $k\in \N$ such that $\frac{1}{2^{k+1}}<t\leq\frac{1}{2^{k}}.$
Then, one can easily see that
\Eq{kk}{
  d_\Z(t)=t,\qquad d_\Z(2t)=2t,\quad\dots,\quad d_\Z(2^{k-1}t)=2^{k-1}t,\qquad
   d_\Z(2^kt)=1-2^kt.
}
On the other hand, by the well-known identity $\sum_{j=0}^{k-1} 2^j=2^k-1$, we have
\Eq{*}
{
(1-t)u=tu+2tu+\cdots+2^{k-1}tu+(1-2^kt)u.
}
Then, by the subadditivity of $\varphi$, and by $t\leq \frac{1}{2^{k}}<\frac{1}{2^{k-1}}<\cdots<\frac{1}{2},$
it follows that
\Eq{*}
{
t\varphi((1-t)u) &\leq t\varphi(tu)+t\varphi(2tu)+\cdots+t\varphi(2^{k-1}tu)+t\varphi((1-2^kt)u)\\
   &\leq t\varphi(tu)+\frac{\varphi(2tu)}{2}+\cdots
           +\frac{\varphi(2^{k-1}tu)}{2^{k-1}}+\frac{\varphi((1-2^kt)u)}{2^k}.
}
Adding $(1-t)\varphi(tu)$ to the previous inequality and using \eq{kk}, we get
\Eq{*}
{
\tau_\varphi(t,u):=t\varphi((1-t)u)+(1-t)\varphi(tu)\leq& \varphi(tu)+\frac{\varphi(2tu)}{2}+\cdots+\frac{\varphi(2^{k-1}tu)}{2^{k-1}}+\frac{\varphi((1-2^kt)u)}{2^k}\\
  =& \sum_{j=0}^k \frac{\varphi(d_\Z(2^jt)u)}{2^j}\leq \T_\varphi(t,u).
}
Which completes the proof of \eq{7a}.
\end{proof}

\Prp{8}{Let $\varphi:D^+\to \R_+$ be nondecreasing with $\varphi(s)>0$ for $s>0$
and assume that
\Eq{*}{
  \gamma_\varphi:=\sup_{0<s\in\frac{1}{2}D^+}\frac{\varphi(2s)}{\varphi(s)}<2.
}
Then, for all $(t,u)\in\R\times D^+$,
\Eq{8a}
{
\T_\varphi(t,u)\leq \frac{2}{2-\gamma_\varphi}\tau_\varphi(t,u)
}
holds.
}

\begin{proof}
To prove \eq{8a}, we fix an arbitrary element $u\in D^+$.
By \rem{1}, the function $\T_\varphi(\cdot,u)$ is the fixed point of the operator
\Eq{*}
{
(F_\varphi f)(t)=\frac{1}{2}f(2t)+\varphi(d_\Z(t)u).
}
Define the function $g:\R\to\R$ by $g(t):=  \frac{2}{2-\gamma_\varphi}\tau_\varphi(t,u)$.
In view of \thm{TT}, in order to prove inequality \eq{8a}, it is enough to show that
\Eq{8b}
{
(F_\varphi g)(t)\leq g(t)\qquad(t\in\R).
}
Since $g$ is periodic by $1$ and symmetric with respect to $t=1/2$,
it suffices to prove that \eq{8b} is satisfied on $\big[0,\frac12\big].$
Trivially, $\gamma_\varphi\geq1$, hence the inequality \eq{8b} is obvious
for $t=0$ or for $u=0$. Thus, we may assume that $u>0$ and $0<t\leq \frac12.$
By the definition of the constant $\gamma_\varphi$, we have that
\Eq{8c}
{
 \varphi(tu)\Big(1-\frac{\gamma_\varphi}{2}\Big) \leq \varphi(tu)-\frac{\varphi(2tu)}{2}.
}
Since $t\leq 2t$ and $1-2t\leq 1-t$ and $\varphi$ is nondecreasing we also have that
\Eq{*}
{
0\leq t(\varphi((1-t)u)-\varphi((1-2t)u))+t\varphi(2tu)-t\varphi(tu).
}
Adding $\varphi(tu)-\dfrac{\varphi(2tu)}{2}$ to the previous inequality and also
using \eq{8c}, we obtain
\Eq{*}
{
\varphi(tu)\Big(1-\frac{\gamma_\varphi}{2}\Big) &\leq
\varphi(tu)-\frac{\varphi(2tu)}{2}\\&\leq t\big(\varphi((1-t)u)-\varphi((1-2t)u)\big)+(1-t)\varphi(tu)-\Big(\frac12-t\Big)\varphi(2tu).
}
Rearranging this inequality, we finally obtain that
\Eq{*}
{
\frac{1}{2-\gamma_\varphi}\big(2t\varphi\big((1-2t)u\big)+(1-2t)\varphi\big(2tu\big)\big)+\varphi(tu)
  \leq \frac{2}{2-\gamma_\varphi}\big(t\varphi\big((1-t)u\big)+(1-t)\varphi\big(tu\big)\big),
}
which means that \eq{8b} is satisfied for all $0<t\leq \frac12.$
\end{proof}

Let $\mu$ be a nonnegative finite Borel measure on $[0,1]$ and let $\supp \mu$ denote
the support of $\mu.$

\Lem{10}{Let $\mu$ be a nonnegative and nonzero finite Borel measure on $[0,1]$
and let $\chi:]0,\infty[\to \R_+$ be defined by
\Eq{*}
{
\chi(s)=\frac{\int_{[0,1]}(2s)^pd\mu(p)}{\int_{[0,1]}s^pd\mu(p)}.
}
Then $\chi$ is nondecreasing on $]0,\infty[$ and
\Eq{lim}{
  \lim_{s\to\infty}\chi(s)=2^{p_0},
}
where $p_0:=\sup(\supp\mu)$.}

\begin{proof}The function $x\mapsto 2^x$ is strictly increasing, hence, for
$p,q\in\R$, we have $ (2^p-2^q)(p-q)\geq0$. It suffices to show that $\chi'\geq0$.
For $s>0,$ we obtain
\Eq{*}
{
\chi'(s)&=\frac{\int_{[0,1]}2^pps^{p-1}d\mu(p) \cdot \int_{[0,1]}s^{p}d\mu(p)
         - \int_{[0,1]}2^ps^{p}d\mu(p) \cdot \int_{[0,1]}ps^{p-1}d\mu(p)}
               {\big(\int_{[0,1]}s^{p}d\mu(p)\big)^2}\\
      &=\frac{\int_{[0,1]}2^pps^{p-1}d\mu(p) \cdot \int_{[0,1]}s^{q}d\mu(q)
         + \int_{[0,1]}2^qqs^{q-1}d\mu(q) \cdot \int_{[0,1]}s^{p}d\mu(p)}
               {2\big(\int_{[0,1]}s^{p}d\mu(p)\big)^2}\\
      &\quad- \frac{\int_{[0,1]}2^ps^{p}d\mu(p) \cdot \int_{[0,1]}qs^{q-1}d\mu(q)
         + \int_{[0,1]}2^qs^{q}d\mu(q) \cdot \int_{[0,1]}ps^{p-1}d\mu(p)}
               {2\big(\int_{[0,1]}s^{p}d\mu(p)\big)^2}\\
        &=\frac{\int_{[0,1]}\int_{[0,1]}(2^p-2^q)(p-q)s^{p+q-1} d\mu(p)d\mu(q)}
               {2\big(\int_{[0,1]}s^{p}d\mu(p)\big)^2}\geq 0,
}
which proves that $\chi$ is nondecreasing.

Using $\supp\mu\subseteq[0,p_0]$, for $s>0,$ we obtain
\Eq{*}
{
\int_{[0,1]}s^pd\mu(p)=\int_{[0,1]}2^p\Big(\frac{s}{2}\Big)^pd\mu(p)\leq 2^{p_0}\int_{[0,1]}\Big(\frac{s}{2}\Big)^pd\mu(p),
}
which proves that $\chi(s)\leq 2^{p_0}$, and hence, $\lim_{s\to\infty}\chi(s)\leq2^{p_0}$.

To show that in \eq{lim} the equality is valid,
assume that $\lim_{s\to\infty}\chi(s)<2^{p_0}$. Choose $q<q_0<p_0$
so that $\lim_{s\to\infty}\chi(s)\leq 2^{q}$. Then, for all $s>0$,
\Eq{*}{
  \int_{[0,1]}(2s)^pd\mu(p)\leq2^q\int_{[0,1]}s^pd\mu(p),
}
i.e., for all $s\geq1$,
\Eq{*}{
  0&\leq\int_{[0,1]}(2^q-2^p)s^pd\mu(p)\\
   &=\int_{[0,q[}(2^q-2^p)s^pd\mu(p)
   +\int_{[q,q_0[}(2^q-2^p)s^pd\mu(p)+\int_{[q_0,1]}(2^q-2^p)s^pd\mu(p)\\
   &\leq\int_{[0,q[}(2^q-2^p)s^pd\mu(p)+\int_{[q_0,1]}(2^q-2^p)s^pd\mu(p)\\
   &\leq\int_{[0,q[}(2^q-2^p)s^pd\mu(p)+\int_{[q_0,1]}(2^q-2^p)s^{q_0}d\mu(p).
}
Therefore, for $s\geq1$,
\Eq{*}{
  0\leq\int_{[0,q[}(2^q-2^p)s^{p-q_0}d\mu(p)+\int_{[q_0,1]}(2^q-2^p)d\mu(p).
}
The first integrand converges uniformly to $0$ on $[0,q[$ as $s\to\infty$.
Thus, by taking the limit $s\to\infty$, we get
\Eq{xx}{
  0\leq\int_{[q_0,1]}(2^q-2^p)d\mu(p).
}
On the other hand, the inequality $q_0<p_0=\sup(\supp\mu)$
implies $\mu([q_0,1])>0$ and, obviously, $2^q-2^p<0$ for $p\in[q_0,1]$. Hence
the right hand side of \eq{xx} is negative. The contradiction so obtained proves
\eq{lim}.
\end{proof}

\Prp{9}{Let $\mu$ be a nonnegative and nonzero finite Borel measure on $[0,1]$.
Denote $\alpha:=\sup D^+$ and $p_0:=\sup(\supp\mu)$ and define $\varphi:D^+\to \R_+$ by
\Eq{*}
{
\varphi(s):=\int_{[0,1]}s^pd\mu(p) \quad \mbox{for all} \quad s\in D^+.
}
Then $\varphi$ is subadditive and nondecreasing, furthermore,
\Eq{gp}
{
\gamma_{\varphi}=\left\{
\begin{array}{lcl}
\dfrac{\int_{[0,1]}\alpha^pd\mu(p)}{\int_{[0,1]}(\alpha/2)^pd\mu(p)},
    &\text{if} & \alpha<\infty,\\[6mm]
2^{p_0},
    &\text{if} & \alpha=\infty
\end{array}\right.
}
and $\gamma_\varphi<2$ if either $\alpha<\infty$ and $\mu$ is not concentrated at the
singleton $\{1\}$ or $p_0<1$.
In addition, for all $t\in[0,1]$ and $u\in D^+$,
\Eq{9A}
{
\int_{[0,1]}\big[t(1-t)^p+(1-t)t^p\big]u^pd\mu(p)
 \leq \int_{[0,1]}T_p(t)u^pd\mu(p)
}
and, provided that $\gamma_\varphi<2$,
\Eq{9B}
{
 \int_{[0,1]}T_p(t)u^pd\mu(p)
 \leq \frac{2}{2-\gamma_\varphi}\int_{[0,1]}\big[t(1-t)^p+(1-t)t^p\big]u^pd\mu(p).
}}

\begin{proof}It can be easily seen that $\varphi$ is nondecreasing. The subadditivity
is a consequence of \lem{4}.

Let $\alpha<\infty.$ Then, by \lem{10}, the map
$s\mapsto\dfrac{\varphi(2s)}{\varphi(s)}= \dfrac{\int_{[0,1]}(2s^p)d\mu(p)}{\int_{[0,1]}s^pd\mu(p)}=\chi(s)$
is nondecreasing on $\frac12D^+,$ so it attains its supremum at $\alpha/2.$ Thus, in this case,
$\gamma_\varphi=\dfrac{\int_{[0,1]}\alpha^pd\mu(p)}{\int_{[0,1]}(\alpha/2)^pd\mu(p)}.$
To prove that $\gamma_\varphi<2$, we use the inequality $2^p<2$ for $p\in[0,1]$ to obtain:
\Eq{*}
{
\int_{[0,1]}\alpha^pd\mu(p)
   =\int_{[0,1]}2^p\Big(\frac{\alpha}{2}\Big)^pd\mu(p)
   <2\int_{[0,1]}\Big(\frac{\alpha}{2}\Big)^pd\mu(p).
}
In the case $\alpha=\infty$, by \lem{10}, we have that
$\gamma_\varphi=\lim_{s\to\infty}\chi(s)=2^{p_0}$.
Obviously, $\gamma_{\varphi}<2$ if $p_0<1$.

The inequalities \eq{9A} and \eq{9B} are immediate consequences of \prp{7} and \prp{8}, respectively.
\end{proof}

In the case when the measure $\mu$ is concentrated at a singleton $\{p\}$, \prp{9}
simplifies to the following result.

\Cor{13}{Let $0\leq p\leq1$ be an arbitrary constant. Then, for all $t\in[0,1]$,
\Eq{*}
{
t(1-t)^p+(1-t)t^p\leq T_p(t)
}
and, provided that $p<1$,
\Eq{*}
{
 T_p(t)\leq \frac{2}{2-2^p}\big(t(1-t)^p+(1-t)t^p\big).
}
}

The proof of the next theorem is analogous to that in \cite{HazPal05}.

\Thm{10}{Let $\varphi:D^+\to\R_+$ be nondecreasing. If $f : D\to \R$ is $\varphi$-midconvex and locally
bounded from above at a point of $D$, then $f$ is locally bounded from above on $D.$}

The following theorem generalizes the analogous result of the paper \cite{HazPal05}
obtained for $(\epsilon,p)$-convexity. A similar result was also established by Tabor and Tabor
\cite{TabTab09b}, \cite{TabTab09a}.

\Thm{H7}{Let $f:D\to \R$ be locally bounded from above at a point of $D$ and
let $\varphi:\frac12 D^+\to\R_+$ be nondecreasing. Then $f$ is $\varphi$-midconvex on $D$,
i.e., \eq{1b} holds for all $x,y\in D$ if and only if
\Eq{H7a}
{
f(tx+(1-t)y)\leq tf(x)+(1-t)f(y)+\T_\varphi(t,\|x-y\|)
}
for all $x,y\in D$ and $t\in [0,1].$}

\begin{proof} Assume that $f$ is $\varphi$-midconvex on $D$ and locally bounded from above
at a point of $D$. From \thm{10}, it follows that $f$ is locally bounded from above
at each point of $D.$
Thus $f$ is bounded from above on each compact subset of $D,$ in particular, for each
fixed $x,y\in D$, $f$ is bounded from above on $[x,y]=\{tx+(1-t)y\mid t\in[0, 1]\}$.
Denote by $K_{x,y}$ a finite upper bound of the function
\Eq{11b}
{
t\mapsto f(tx +(1-t)y)-tf(x)-(1-t)f(y) \qquad (t \in [0, 1]).
}
We are going to show, by induction on $n$, that
\Eq{11c}
{
f(tx+(1-t)y) \leq tf(x)+(1-t)f(y)+\frac{K_{x,y}}{2^n}
             +\sum_{j=0}^{n-1}\frac{\varphi\big(d_\Z(2^{j}t)\|x-y\|\big)}{2^{j}}
}
for all $x, y \in D$ and $t \in [0, 1].$
For $n = 0,$ the statement follows from the definition of $K_{x,y}$ (with the convention
that the summation for $j=0$ to $(-1)$ is equal to zero).

Now assume that \eq{11c} is true for some $n \in \N.$ Assume that $t\in [0,1/2].$
Then, due to the $\varphi$-midconvexity of $f,$ we get
\Eq{*}
{
f(tx+(1-t)y) = f \Big(\frac{y+(2tx+(1-2t)y)}{2}\Big)
\leq \frac{f(y)+f(2tx+(1-2t)y)}{2}+\varphi(t\|x-y\|).
}
On the other hand, by $\eq{11c}$, we get that
\Eq{*}
{
f(2tx+(1-2t)y)\leq 2tf(x)+(1-2t)f(y)+\frac{K_{x,y}}{2^n}
   +\sum_{j=0}^{n-1}\frac{\varphi\big(d_\Z(2^{j+1}t)\|x-y\|\big)}{2^{j}}.
}
Combining these two inequalities, we obtain
\Eq{*}
{
f(tx+(1-t)y)&\leq tf(x)+ (1-t)f(y)+\frac{1}{2}\bigg(\frac{K_{x,y}}{2^n}
          +\sum_{j=0}^{n-1}\frac{\varphi\big(d_\Z(2^{j+1}t)\|x-y\|\big)}{2^{j}}\bigg)+\varphi(t\|x-y\|)\\
            &= tf(x) + (1-t)f(y)+\frac{K_{x,y}}{2^{n+1}}
              +\sum_{j=0}^{n}\frac{\varphi\big(d_\Z(2^{j}t)\|x-y\|\big)}{2^{j}}.
}
In the case $t\in[1/2,1]$, the proof is similar.
Thus, \eq{11c} is proved for all $n\in\N$.

Finally, taking the limit $n\to \infty$ in \eq{11c}, we get the desired inequality \eq{H7a}.

To see that \eq{H7a} implies the $\varphi$-midconvexity of $f$, substitute $t=1/2$ into \eq{H7a}
and use the easy-to-see identity $\T_\varphi\big(\frac12,u\big)=\varphi(\frac{|u|}{2})$ ($u\in\R$).
\end{proof}

The optimality of the error term in \eq{H7a} and the appropriate convexity properties of
$\T_\varphi$ have recently been obtained in \cite{MakPal10b}.

\Thm{8}{Let $\varphi:D^+\to \R_+$ be nondecreasing with $\varphi(s)>0$ for $s>0$
and assume that $\gamma_\varphi:=\sup_{0<s\in\frac12D^+}\frac{\varphi(2s)}{\varphi(s)}<2.$
If $f:D\to \R$ is locally bounded from above a point of $D$ and it is also
$\varphi$-midconvex, then $f$ is $\big(\frac{2}{2-\gamma_\varphi}\cdot\varphi\big)$-convex on $D$.}

\begin{proof}
By \prp{8} and by \thm{H7}, the proof of this theorem is evident.
\end{proof}

\Cor{9}{Let $\mu$ be a nonnegative and nonzero finite Borel measure on $[0,1]$.
Denote $\alpha:=\sup D^+$ and $p_0:=\sup(\supp\mu)$
and assume that either $\alpha<\infty$ and $\mu$ is not concentrated at the
singleton $\{1\}$ or $p_0<1$. Define $\varphi:D^+\to \R_+$ by
\Eq{*}
{
\varphi(s):=\int_{[0,1]}s^pd\mu(p) \quad \mbox{for all} \quad s\in D^+.
}
If $f:D\to \R$ is locally bounded from above a point of $D$ and it is also $\varphi$-midconvex,
then $f$ is $\big(\frac{2}{2-\gamma_{\varphi}}\cdot \varphi\big)$-convex on $D$, where $\gamma_\varphi$
is given by \eq{gp}.}

\Cor{10}{Let $0\leq p<1$ and $\epsilon\geq 0$ be arbitrary constants. 
If $f:D\to \R$ is locally bounded from above a point of $D$ and it is also $(\epsilon,p)$-midconvex, 
then $f$ is $\big(\frac{2\epsilon}{2-2^p},p\big)$-convex on $D$.}


\end{document}